\documentclass[12pt]{article}
\usepackage{latexsym}
\usepackage{amsmath}
\usepackage{amssymb}
\usepackage{amsthm}
\usepackage{stmaryrd}
\usepackage{epsfig}
\usepackage{graphicx}
\usepackage{enumerate}
\usepackage{subfigure}
\usepackage{url}
\usepackage{bbm}

\theoremstyle{plain}
\newtheorem{thm}{Theorem}[section]
\newtheorem{lem}[thm]{Lemma}

\theoremstyle{definition}

\theoremstyle{remark}

\newcommand{\be}{\begin{enumerate}}
\newcommand{\ee}{\end{enumerate}}
\newcommand{\smn}{\setminus}

\advance \topmargin by -\headheight
\advance \topmargin by -\headsep
\evensidemargin \oddsidemargin
\usepackage[pdftex,
    pdftitle={A Linear Programming Relaxation and a Heuristic for the Restless Bandit Problem with General Switching Costs},
    pdfauthor={Jerome Le Ny},
    colorlinks,linkcolor={blue},citecolor={blue},urlcolor={red},
]
{hyperref}

\title{\LARGE \bf
A Linear Programming Relaxation and a Heuristic for the Restless Bandits Problem with General Switching Costs
}

\author{ \parbox{4 in}{\centering Jerome Le Ny, Munther Dahleh and Eric Feron}
         \thanks{This work was supported by Air Force - DARPA - MURI award 009628-001-03-132 and Navy ONR award N00014-03-1-0171. J. Le Ny and M. Dahleh are with the Laboratory for Information and Decision Systems, Massachusetts Institute of Technology, Cambridge, MA 02139-4307, USA. E. Feron is with the School of Aerospace Engineering, Georgia Tech, Atlanta, GA 30332, USA. {\tt\small jleny@mit.edu, eric.feron@aerospace.gatech.edu} }\\
}

\begin{document}

\maketitle
\thispagestyle{empty}
\pagestyle{empty}

\begin{abstract}

We extend a relaxation technique due to Bertsimas and Ni{\~n}o-Mora for the restless bandit problem to the case where arbitrary costs penalize switching between the bandits. We also construct a one-step lookahead policy using the solution of the relaxation. Computational experiments and a bound for approximate dynamic programming provide some empirical support for the heuristic.  

\end{abstract}





\section{Introduction} \label{intro} 

We study the restless bandit problem (RBP) with general switching costs between the bandits, which could represent travel distances for example. This problem is an intractable extension of the multi-armed bandit problem (MABP), which can be described as follows. There are $N$ projects, of which only one can be worked on at any time period. Project $i$ is characterized at (discrete) time $t$ by its state $x_i(t)$, which belongs to a finite state space $S_i$. If project $i$ is worked on at time $t$, one receives a reward $\alpha^t r(x_i(t))$, where $\alpha \in (0,1)$ is a discount factor. The state $x_i(t)$ then evolves to a new state according to given transition probabilities. The states of all idle projects are unaffected. The goal is to find a policy which decides at each time period which project to work on in order to maximize the expected sum of the discounted rewards over an infinite horizon. The MABP problem was first solved by Gittins \cite{gittinsJones}. He showed that it is possible to define separately for each project an index which is a function of the project state only, and that the optimal policy operates at each period the project with the greatest current index. Moreover, these indices can be calculated efficiently, as shown for example in \cite{varaiya-extensions}.

Whittle \cite{whittle-restless} proposed an interesting modification of the model, called the restless bandit problem (RBP), which extends significantly the range of applications. In the RBP, one can activate several projects at each time period, and the projects that are not activated continue to evolve, possibly using different transition probabilities. Finding an optimal policy efficiently for the RBP is unlikely to be possible however, since the problem is PSPACE-hard \cite{papadimitriou-complexityQueuing}, even in restricted cases. Nevertheless, Whittle proposed an index policy for the RBP which performs well in practice. 

Another extension of the MABP concerns the addition of costs for changing the currently active project. This problem, which we call the multi-armed bandit problem with switching costs (MABPSC), is of great interest to various applications, as discussed by \cite{jun-surveypaper}, \cite{jun05-investments}, \cite{vanOyen-indexPolicies}, \cite{koole-server}, in order to model for example set-up and tear-down costs in queuing networks, transition costs in a job search problem or transaction fees in a portfolio optimization problem. It is easy to see that the MABPSC is NP-hard, since the HAMILTON CYCLE problem is a special case of it \cite{leny06_RBSC_ACC}. The MABPSC has been studied in particular by Asawa and Teneketzis \cite{asawa-penalties}, and very recently by Glazebrook et al. \cite{Glazebrook06_machine_index} and Ni{\~n}o-Mora \cite{ninoMora-MABPSC}. These authors are concerned with the case where the switching costs have a separable form $c_{ij}=c_i+c_j$, preserving the separable structure from the MABP, and design approximate index policies.

Our work was motivated by an optimal aerial surveillance problem, where switching costs correspond to travel distances between inspection sites. Hence, the assumption on the separable form of the switching costs does not hold. This introduces additional coupling between the projects, and it is not clear then how to design index policies. Moreover, the sites continue to evolve while not visited, and thus we are led to consider the restless bandit problem with switching costs (RBPSC). 

We adopt a computational approach to the RBPSC. We impose no restriction on the switching costs, not even the triangle inequality. In Section \ref{MDP}, we formulate the problem as a Markov decision process (MDP), using the state-action frequency approach \cite{derman-LP}. This yields a linear program, which we relax in section \ref{relaxation} by following an idea that Bertsimas and Ni{\~n}o-Mora developped for the RBP \cite{bertsimas-RB}, optimizing over a restricted set of marginals of the occupation measure. The coupling introduced by the switching costs makes this relaxation significantly more challenging to develop than in the classical case, and the first contribution of the paper is to present valid constraints on the marginals improving the quality of the relaxation. This relaxation provides an efficiently computable bound on the achievable performance. Section \ref{algorithms} describes how the relaxation can also be used to motivate a heuristic policy. This heuristic is based on approximate dynamic programming (ADP) techniques, but we also show how to recover it from the linear programming theory used by Bertsimas and Ni{\~n}o-Mora to design their primal-dual heuristic for the RBP. Section \ref{experiments} presents numerical experiments comparing the heuristic to the performance bound. 

The advantage of using the approximate dynamic programming point of view is that a recently developed performance bound provides additional support for our heuristic. However, we do not consider in this paper the development of policies with an priori performance bound. Few results exist in the literature concerning such bounds. As remarked by Guha et al. \cite{Guha07_RBApprox}, even the standard RBP is PSPACE-Hard to approximate to any non-trivial factor, unless some assumptions are made on the reward functions.

\section{Exact Formulation of the RBSC Problem} \label{MDP}

We formulate the RBPSC using the linear programming approach to Markov decision processes \cite{derman-LP}, \cite{puterman-MDP}. $N$ projects are distributed in space at $N$ sites, and $M \leq N$ servers can be allocated to $M$ different projects at each time period $t=1,2,\ldots$. In the following, we use the terms project and site interchangeably; likewise, agent and server have the same meaning. At each time period, each server must occupy one site, and different servers must occupy distinct sites. We say that a site is active at time $t$ if it is visited by a server, and is passive otherwise. If a server travels from site $k$ to site $l$, we incur a cost $c_{kl}$. Each site can be in one of a finite number of states $x_n \in S_n$, for $n=1,\ldots,N$, and we denote the Cartesian product of the individual state spaces by $\mathcal{S}=S_1 \times \ldots \times S_N$. If site $n$ in state $x_n$ is visited, a reward $r^1_n(x_n)$ is earned, and its state changes to $y_n$ according to the transition probabilities $p^1_{x_n y_n}$. If the site is not visited, then a reward (potentially negative) $r^0_n(x_n)$ is earned for that site and its state changes according to the transition probabilities $p^0_{x_n y_n}$. \emph{We assume that all sites change their states independently of each other}. 

Let us denote the set $\{1,\ldots,N\}$ by $[N]$. We consider that when no agent is present at a given site, there is a fictitious agent called passive agent at that site. We also call the real agents active agents, since they collect active rewards. The transition of a passive agent between sites does not involve any switching cost, and when a passive agent is present at a site, the passive reward is earned. Therefore, we have a total of $N$ agents including both the real and passive agents, and we can describe the positions of all agents by a vector $\mathbf{s}=(s_1,\ldots,s_N)$, which corresponds to a \emph{permutation} of [N]. 
We denote the set of these permutation vectors by $\Pi_{[N]}$, and we let the $M$ first components correspond to the active agents. For example, with $M=2$ and $N=4$, the vector $(s_1=2,s_2=3,s_3=1,s_4=4) \in \Pi_{[4]}$ means that agent $1$ is in site $2$, agent $2$ in site $3$ and sites $1$ and $4$ are passive. 

For an agent $i \in [N]$, we refer to the set of the other agents by $-i$. If we fix $s_i \in [N]$ for some $1 \leq i \leq N$, then we write $\mathbf{s}_{-i}$ to denote the vector $(s_1,\ldots,s_{i-1},s_{i+1},\ldots,s_N)$, and $\Pi_{[N]-s_i}$ to denote the permutations of the set $[N]-\{s_i\}$. In particular, we write $\sum_{\mathbf{s}_{-i} \in \Pi_{[N]-s_i}}$ to denote the sum over all permutations of the positions of the agents $-i$, over the set of sites not occupied by agent $i$. We also write $\mathcal{S}_{-i}$ to denote the cartesian product $S_1 \times \ldots S_{i-1} \times S_{i+1} \times \ldots \times S_{N}$.

The state of the system at time $t$ can be described by the state of each site and the positions $\mathbf{s} \in \Pi_{[N]}$ of the servers, including the passive ones. With this state description, we are able to handle any number $M \leq N$ of agents as a parameter within the same framework. We denote the complete state by $(x_1,\ldots,x_N;s_1,\ldots,s_N):=(\mathbf{x};\mathbf{s})$. We can choose which sites are to be visited next, i.e., an action $\mathbf{a}$ belongs to the set $\Pi_{[N]}$ and corresponds to the assignment of the agents, including the passive ones, to the sites for the next time period. Once the sites to be visited are chosen, there are costs $c_{s_i a_i}$ for moving the active agent $i$ from site $s_i$ to site $a_i$, including possibly a nonzero cost for staying at the same site. The immediate reward earned is
\[
R((\mathbf{x};\mathbf{s}),\mathbf{a}):=\sum_{i=1}^M (r^1_{a_i}(x_{a_i})-c_{s_i a_i})+\sum_{i=M+1}^N r^0_{a_i}(x_{a_i}).
\] 
We are given a distribution $\nu$ on the initial state of the system, and \emph{we will assume a product form} 
\begin{equation} \label{eq: initial distribution - product}
\nu(\mathbf{x};\mathbf{s})=\prod_{i=1}^N \nu_i(x_i) \delta_{d_i}(s_i),
\end{equation}
i.e., the initial states of the sites are independent random variables and server $i$ leaves initially from site $d_i$, with $\mathbf{d} \in \Pi_{[N]}$. 


The transition probability matrix has a particular structure, since the sites evolve independently and the transitions of the agents are deterministic. Let us write its elements $\mathcal{P}_{(\mathbf{x'};\mathbf{s'})\mathbf{a}(\mathbf{x};\mathbf{s})}=\mathcal{P}_{\mathbf{x} \, \mathbf{a} \, \mathbf{x'}} \prod_{i=1}^N \delta_{s_i}(a_i)$, where
\[
\mathcal{P}_{\mathbf{x} \, \mathbf{a} \, \mathbf{x'}} := \prod_{i=1}^M p^1_{x_{a_i}x'_{a_i}} \prod_{i=M+1}^N p^0_{x_{a_{i}}x'_{a_{i}}}.
\]

The optimal infinite horizon discounted reward, multiplied by $(1-\alpha)$, is the optimal value of the following linear program (LP) \cite{derman-LP}
\begin{align} 
& \text{maximize } \sum_{\mathbf{s} \in \Pi_{[N]}} \sum_{\mathbf{a} \in \Pi_{[N]}} \sum_{\mathbf{x} \in \mathcal{S}} R((\mathbf{x};\mathbf{s}),\mathbf{a}) \, \rho_{(\mathbf{x};\mathbf{s}),\mathbf{a}}  \label{RBSC Multi-agent} \\ 
& \text{subject to } \nonumber \\
& \sum_{\mathbf{s}',\mathbf{a} \in \Pi_{[N]}} \sum_{\mathbf{x'} \in \mathcal{S}}  \rho_{(\mathbf{x}';\mathbf{s}'),\mathbf{a}} [\delta_{(\mathbf{x};\mathbf{s})}(\mathbf{x'};\mathbf{s'}) - \alpha \mathcal{P}_{(\mathbf{x'};\mathbf{s'})\mathbf{a}(\mathbf{x};\mathbf{s})}] \nonumber \\
& =(1-\alpha) \nu(\mathbf{x};\mathbf{s}), \quad \forall \; (\mathbf{x},\mathbf{s}) \in \mathcal{S}\times \Pi_{[N]} \label{RBSC Multi-agent constraints} \\
& \rho_{(\mathbf{x};\mathbf{s}),\mathbf{a}} \geq 0, \quad \forall \; ((\mathbf{x};\mathbf{s}),\mathbf{a}) \in \mathcal{S}\times \Pi_{[N]}^2. \nonumber
\end{align}

The variables $\{\rho_{(\mathbf{x};\mathbf{s}),\mathbf{a}}\}$ of the LP, called state action frequencies or occupation measure, form a probability measure on the space of state-action pairs and an optimal policy can be recovered from an optimal solution for the LP. The formulation above is of little computational interest however since the number of variables and constraints is of the order of $|\mathcal{S}| \times (N!)^2$, that is, exponential in the size of the input. 

We can obtain the linear program dual to (\ref{RBSC Multi-agent}) by constructing it directly, or starting from Bellman's equation and using standard dynamic programming arguments \cite[vol. 2, p. 53]{bertsekas-dynamicProgramming}. The decision variables $\{ \lambda_{\mathbf{x},\mathbf{s}} \}_{\mathbf{x},\mathbf{s}}$ of the dual correspond to the reward-to-go vector. We get
\begin{align} 
& \text{minimize } (1-\alpha) \sum_{ \mathbf{x},\mathbf{s} } \lambda_{\mathbf{x},\mathbf{s}} \nu(\mathbf{x};\mathbf{s})  \label{eq: dual exact RBSC} \\
& \text{s.t. } \lambda_{\mathbf{x},\mathbf{s}} - \alpha \sum_{ \mathbf{\tilde x}  \in \mathcal{S} } \mathcal{P}_{\mathbf{x} \, \mathbf{a} \, \mathbf{\tilde x}} \,  \lambda_{\mathbf{\tilde x},\mathbf{a}} \geq R((\mathbf{x};\mathbf{s}),\mathbf{a}) \; , \nonumber \\
&\forall \mathbf{x} \in \mathcal{S}, \forall (\mathbf{s},\mathbf{a}) \in \Pi_{[N]}^2. \nonumber
\end{align}


\section{LP Relaxation of the RBPSC} \label{relaxation}

We compute a bound on the performance achievable by any assignment policy by relaxing the LP formulation. We start by rewriting the objective function (\ref{RBSC Multi-agent}) and we identify the relevant marginals of the occupation measure:
\begin{align*}
&\sum_{s=1}^N \sum_{a=1}^N \sum_{x_{a} \in S_{a}} \left[ \left( r^1_{a}(x_{a}) - c_{sa} \right) \left( \sum_{i=1}^M \rho^i_{x_{a};s,a} \right) \right. \left. + r^0_{a}(x_{a}) \left( \sum_{i=M+1}^N \rho^i_{x_{a};s,a} \right) \right] 
\end{align*}
where the marginals appearing above are obtained as follows:
\begin{align} 
& \rho^{i}_{(x_{a};s),a}= \sum_{\mathbf{x}_{-a} \in \mathcal{S}_{-a}} \; \sum_{\mathbf{s}_{-i} \in \Pi_{[N]-s}} \; \sum_{\mathbf{a}_{-i} \in \Pi_{[N]-a}} \;  \rho_{(\mathbf{x};\mathbf{s}),\mathbf{a}} \; , \; \nonumber \\
& \forall \; (i,s,a) \in [N]^3, \; x_a \in S_a, \label{variables needed I bis}
\end{align}
and the superscripts refer to the agents. 

Now to express the constraints, we will also need the following variables:
\begin{align} 
& \rho^{i}_{(x_{s};s),a}= \sum_{\mathbf{x}_{-s} \in \mathcal{S}_{-s}} \; \sum_{\mathbf{s}_{-i} \in \Pi_{[N]-s}} \; \sum_{\mathbf{a}_{-i} \in \Pi_{[N]-a}} \;  \rho_{(\mathbf{x};\mathbf{s}),\mathbf{a}} \; , \; \nonumber \\
& \forall \; (i,s,a) \in [N]^3, \; x_s \in S_s. \label{variables needed I}
\end{align}
The variables in (\ref{variables needed I bis}) (respectively (\ref{variables needed I})) can be interpreted as the frequency with which agent $i$ switches from site $s$ to site $a$ and the destination (resp. origin) site is in state $x_a$ (resp. $x_s$). 
Note that this notation is somewhat redundant, since we can write the variables $\rho^{i}_{(x_{j};j),j}$ as in (\ref{variables needed I bis}) or (\ref{variables needed I}).

It is straightforward to see that the constraints (\ref{RBSC Multi-agent constraints}) imply
\begin{align}
& \sum_{a=1}^N \rho^{i}_{(x_{s};s),a}- \alpha \sum_{\tilde x_{s} \in S_{s}} \sum_{s'=1}^N \rho^{i}_{(\tilde x_{s};s'),s} \, p^{{1}\{i \leq M\}}_{\tilde x_{s} x_{s}} \label{st0} \\
& = (1-\alpha) \, \nu_{s}(x_{s}) \delta_{d_i}(s)\; , \; \forall \; (i,s) \in [N]^2, \forall x_{s} \in S_{s}, \nonumber
\end{align}
on the marginals, where ${1} \{\cdot \}$ is the indicator function. However, there are additional relations that must exist because the marginals are obtained from the same original occupation measure. These relations must be found to insure a sufficiently strong relaxation. Another intuitive way to think about this type of constraints is that we enforce sample path constraints only in average \cite{whittle-restless}. First, from the definitions we have immediately:
\begin{align} 
\sum_{x_{s} \in S_{s}} \rho^{i}_{(x_{s};s),a} = \sum_{x_{a} \in S_{a}} \rho^{i}_{(x_{a};s),a} \, , \, \forall \; (i,s,a) \in [N]^3. \label{compat}
\end{align}

Now for the RBPSC, exactly one agent (active or passive) must be at each site at each period. The frequency with which the agents leave site $j$ in state $x_j$ should be equal to the frequency with which the agents move to site $j$ in state $x_j$. So we expect that the following constraints should hold:
\begin{equation} \label{st1}
\sum_{i=1}^N \sum_{a=1}^N \rho^i_{\tilde x_j,j,a} = \sum_{i=1}^N \sum_{s=1}^N \rho^i_{\tilde x_j,s,j}, \quad \forall j \in [N], \tilde x_j \in S_j.
\end{equation}
We now show that (\ref{st1}) are indeed valid constraints. We use the notation $(\mathbf{x_{-j}},\tilde x_j)$ to mean that the $j^{\text{th}}$ component of the vector is $\tilde x_j$, and similiarly for $(\mathbf{s_{-i}},j)$. 
We have, starting from the definition (\ref{variables needed I})
\begin{align}
\sum_{i=1}^N \sum_{a=1}^N \rho^i_{\tilde x_j,j,a} \nonumber 
&= \sum_{\mathbf{x_{-j}} \in S_{-j}} \; \sum_{\mathbf{a} \in \Pi_{[N]}} \; \sum_{i=1}^N \sum_{\mathbf{s}_{-i} \in \Pi_{[N]-j}} \;  \rho_{(\mathbf{x_{-j}},\tilde x_j);(\mathbf{s_{-i}},j),\mathbf{a}} \nonumber \\
&= \sum_{\mathbf{x_{-j}} \in S_{-j}} \; \sum_{\mathbf{a} \in \Pi_{[N]}} \; \sum_{\mathbf{s} \in \Pi_{[N]}} \;  \rho_{(\mathbf{x_{-j}},\tilde x_j);\mathbf{s},\mathbf{a}}. \label{common quantity} 
\end{align}
The first equality comes from the fact that we count all the permutation vectors $\mathbf{a}$ by varying first the $i^{\text{th}}$ component $a_i$ from $1$ to $N$. The second equality comes from the fact that we count all the permutations $\mathbf{s}$ by varying the position $i$ where the component $s_i$ is equal to $j$ (exactly one of the components of a permutation vector of $\Pi_{[N]}$ has to be $j$). The proof that the right hand side of (\ref{st1}) is also equal to the quantity in (\ref{common quantity}) is identical.

Here are two additional sets of valid constraints:
\begin{align}
 \sum_{s=1}^N \; \sum_{x_s \in S_s} \; \sum_{a \in [N] - \tilde a} \rho^i_{x_s;s,a} & = \sum_{k \in [N] - i} \; \sum^N_{s =1} \; \sum_{x_s \in S_s} \rho^k_{x_s;s;\tilde a} \;, \quad \forall (i,\tilde a) \in [N]^2, \label{st2} \\
 \sum_{a=1}^N \; \sum_{x_a \in S_a} \; \sum_{s \in [N] - \tilde s} \rho^i_{x_a;s,a} & = \sum_{k \in [N] - i} \; \sum^N_{a=1} \; \sum_{x_a \in S_a} \rho^k_{x_a;\tilde s;a} \;, \quad \forall (i,\tilde s) \in [N]^2. \label{st3} 
\end{align}
Intuitively, on the left hand side we have the probability that agent $i$ does not go to site $\tilde a$ (respectively does not leave from site $\tilde s$), which must equal the probability that some other agent $k$ (passive or not) goes to site $\tilde a$ (respectively leaves from site $\tilde s$). Again, these relations can be verified by inspection of (\ref{variables needed I}). Indeed, in (\ref{st2}) for $(i,\tilde a)$ fixed, similarly to (\ref{st1}), we have two equivalent ways of summing the occupation measure over all indices $(\mathbf{x},\mathbf{s},\mathbf{a})$ such that none of the permutation vectors $\mathbf{a}$ with $\tilde a$ in position $i$ appears. On the left hand side of (\ref{st2}), we vary the coefficient $a_i$ in the set $\{1,\ldots,N\} \smn \{\tilde a\}$, whereas on the right hand side, we obtain the same result by forcing the element $\tilde a$ to be in a position different from position $i$. Similarly, in (\ref{st3}), we have two ways of summing over all indices such that none of the permutation vectors $\mathbf{s}$ with $\tilde s$ in position $i$ appears.

Finally we have obtained a relaxation for the RBPSC:

\begin{thm} We can obtain an upper bound on the optimal reward achievable in the RBPSC by solving the following linear program:
\begin{align} 
& \text{maximize }  \nonumber \\
& \sum_{s=1}^N \sum_{a=1}^N \sum_{x_{a} \in S_{a}} \left[ \left( r^1_{a}(x_{a}) - c_{sa} \right) \left( \sum_{i=1}^M \rho^i_{x_{a};s,a} \right) \right. \left. + r^0_{a}(x_{a}) \left( \sum_{i=M+1}^N \rho^i_{x_{a};s,a} \right) \right] \label{relaxation Multi-agent} \\
& \text{subject to } \nonumber \\
& \; (\ref{st0}), (\ref{compat}), (\ref{st1}), (\ref{st2}), \; (\ref{st3}),  \nonumber \\
& \rho^i_{(x_s;s),a} \geq 0 \; , \quad \forall (i,s,a) \in [N]^3, x_s \in S_s \nonumber \\
& \rho^i_{(x_a;s),a} \geq 0 \; , \quad \forall (i,s,a) \in [N]^3, x_a \in S_a. \nonumber
\end{align}
\end{thm}


There are now $\mathcal{O}(N^3 \times \max_i |S_i|)$ variables $\rho^i_{(x_s;s)a}, \rho^i_{(x_a;s)a}$, and constraints in the relaxed linear program, which is polynomial in the size of the input. From the remarks about the complexity of the problem, it is unlikely that a polynomial number of variables will suffice to formulate the RBPSC exactly. However, the addition of the constraints tying the marginals together helps reduce the size of the feasible region spanned by the decision vectors and improve the quality of the relaxation. Computing the optimal value of this linear program can be done in polynomial time, and provides an upper bound on the performance achievable by any policy for the original problem.



\subsection{Dual of the Relaxation} \label{section: relaxation dual}


It will be useful to consider the dual of the LP relaxation obtained in the previous paragraph, which we derive directly from (\ref{relaxation Multi-agent}). This dual program could be obtained by dynamic programming arguments, in the spirit of the original work of Whittle, incorporating the constraints $(\ref{compat}), (\ref{st1}), (\ref{st2}), (\ref{st3})$ using Lagrange multipliers. 
We obtain:
\begin{align}
& \text{minimize } (1-\alpha) \sum_{i=1}^N \sum_{s=1}^N \sum_{x_s \in S_s} \nu_s(x_s) \delta_{d_i}(s) \lambda^i_{s,x_s} \label{eq: relaxation dual multi-agent} \\
& \text{subject to } \nonumber \\
& \lambda^{i}_{s, x_{s}} + \mu^{i}_{s, a} + \kappa_{s, x_{s}} - \sum_{i' \neq i} \zeta^{i'}_{a} + \sum_{a' \neq a} \zeta^{i}_{a'} 
 \geq \, 0 \;, \quad \forall \; (i,s,a) \in [N]^3, s \neq a, \label{eq: constraint xs,s,a} \\
& - \alpha \sum_{\tilde x_{a}} p^{{1}\{i \leq M\}}_{x_{a} \tilde x_{a}} \lambda^{i}_{a, \tilde x_{a}} - \mu^{i}_{s, a} - \kappa_{a, x_{a}} - \sum_{i' \neq i} \xi^{i'}_{s} + \sum_{s' \neq s} \xi^{i}_{s'} \geq r^{{1}\{i \leq M\}}_{a} (x_{a}) - c_{s a} {{1}\{i \leq M\}} \;,  \nonumber \\ 
& \forall \; (i,s,a) \in [N]^3, s \neq a, \label{eq: constraint xa,s,a} \\
& \lambda^{i}_{s, x_{s}} - \alpha \sum_{\tilde x_{s}} p^{{1}\{i \leq M\}}_{x_{s} \tilde x_{s}} \lambda^{i}_{s, \tilde x_{s}} - \sum_{i' \neq i} ( \zeta^{i'}_{s} + \xi^{i'}_{s} ) + \sum_{s' \neq s} ( \zeta^{i}_{s'} + \xi^{i}_{s'} ) \geq r^{{1}\{i \leq M\}}_{s} (x_{s}) - c_{s s} {{1}\{i \leq M\}} \; , \nonumber \\
& \forall i \in [N]. \label{eq: constraint xs,s,s} 
\end{align}

The optimal dual variables $\bar \lambda^i_{s,x_s}$ are Lagrange multipliers corresponding to the constraints (\ref{st0}) and have a natural interpretation in terms of reward-to-go if site $s$ is in state $x_s$ and visited by agent $i$. The optimal dual variables $\bar \mu^i_{s,a}, \bar \kappa_{a,x_a}, \bar \zeta^{i}_{a}, \bar \xi^{i}_{s}$ correspond to the additional constraints (\ref{compat}), (\ref{st1}), (\ref{st2}), and (\ref{st3}) respectively. We can obtain the optimal primal and dual variables simultaneously when solving the relaxation. For $j=s$ or $a$, we also obtain the optimal reduced costs $\bar \gamma^i_{x_j,s,a}$: $\bar \gamma^i_{x_s,s,a}$ are equal to the left hand side of the constraints (\ref{eq: constraint xs,s,a}), whereas $\bar \gamma^i_{x_a,s,a}$ are equal to the difference between the left hand side and the right hand side of (\ref{eq: constraint xa,s,a}), or (\ref{eq: constraint xs,s,s}) if $s=a$. There is one such reduced cost for each variable $\rho^i_{x_j,s,a}$ of the primal, and by complementary slackness, $\bar \rho^i_{x_j,s,a} \bar \gamma^i_{x_j,s,a} = 0$, where $\{ \bar \rho^i_{x_j,s,a} \}$ is the optimal solution of the primal.



\section{A Heuristic for the RBPSC} \label{algorithms}

The relaxation is also useful to actually design assignment policies for the agents. We present here a one-step lookahead policy and its relationship with the primal-dual heuristic of Bertsimas and Nin\~o-Mora, developped for the RBP. 

\subsection{One-Step Lookahead Policy}

Consider the multi-agent system in state $(\mathbf{x};\mathbf{s})$, 
with $\mathbf{s}$ a permutation of $[N]$. Given the interpretation of the dual variables $\lambda^i_{s_i,x_{s_i}}$ in terms of reward-to-go mentioned in section \ref{section: relaxation dual}, it is natural to try to form an approximation $\tilde{J}(\mathbf{x};\mathbf{s})$ of the global reward-to-go in state $(\mathbf{x};\mathbf{s})$ as
\begin{equation} \label{approximate value}
\tilde{J}(x_1,\ldots,x_N;s_1,\ldots,s_N) = \sum_{i=1}^N \bar{\lambda}^i_{s_i,x_{s_i}} ,
\end{equation}
where $\bar{\lambda}^i_{x_{s_i},s_i}$ are the optimal values of the dual variables obtained when solving the LP relaxation. 
The separable form of this approximate cost function is useful to design an easily computable one-step lookahead policy \cite{bertsekas-dynamicProgramming}, as follows. In state $(\mathbf{x};\mathbf{s})$, we obtain the assignment $\tilde{u}(\mathbf{x};\mathbf{s})$ of the agents as
\begin{align}
&\tilde{u}(\mathbf{x};\mathbf{s}) \in \arg \max_{\mathbf{a} \in \Pi_{[N]}} \Bigg \{ R((\mathbf{x};\mathbf{s}),\mathbf{a}) \left. +\alpha \sum_{\mathbf{x'} \in \mathcal{S}} \mathcal{P}_{\mathbf{x} \, \mathbf{a} \, \mathbf{x'}} \, \tilde{J}({\mathbf{x'};\mathbf{a}}) \right\}. 
\end{align}
In this computation, we replaced the true optimal cost function, which would provide an optimal policy, by the approximation $\tilde J$. Using (\ref{approximate value}), we can rewrite the maximization above as
\begin{align} \label{eq: linear assignment}
\max_{\mathbf{a} \in \Pi_{[N]}} \; \sum_{i=1}^N m_{i,a_i},
\end{align}
with 
\begin{align*}
m_{i,a_i} = r^{{1}\{i \leq M\}}_{a_i}(x_{a_i}) - c_{s_i a_i} {{1}\{i \leq M\}} + \alpha \sum_{x'_{a_i} \in S_{a_i}}  \bar{\lambda}^i_{a_i,x'_{a_i}} \, p^{{1}\{i \leq M\}}_{x_{a_i} x'_{a_i}}.
\end{align*}


Assuming that the optimal dual variables have been stored in memory, the evaluation of the terms $m_{i,a_i}$, for all $(i,a_i)$, takes a time $\mathcal{O}(N^2 \max_i{|S_i|})$. The maximization (\ref{eq: linear assignment}) is then a linear assignment problem, which can be solved by linear programming or in time $\mathcal{O}(N^3)$ by the Hungarian method \cite{schrijver-book}. Thus, the assignment can be computed at each time step in time $\mathcal{O}(N^2 \max_i{|S_i|}+N^3)$ by a centralized controller.

\subsection{Equivalence with the Primal-Dual Heuristic}

Recall from paragraph \ref{section: relaxation dual} that, when solving the linear programming relaxation, we can obtain the optimal primal variables $\{ \bar \rho^i_{x_j,s,a} \}$,  the dual variables $\{ \bar \lambda^i_{s,x_s},\bar \mu^i_{s,a}, \bar \kappa_{a,x_a}, \bar \zeta^{i}_{a}, \bar \xi^{i}_{s} \}$, and  the reduced costs $\{ \bar \gamma^i_{x_j,s,a} \}$. These reduced costs are nonnegative. Bertsimas and Nin\~o-Mora \cite{bertsimas-RB} motivated their primal-dual heuristic for the RBP using the following well-known interpretation of the reduced costs: \emph{starting from an optimal solution, $\bar \gamma^i_{(x_j;s),a}$ is the rate of decrease in the objective value of the primal linear program (\ref{relaxation Multi-agent}) per unit increase in the value of the variable $\rho^i_{(x_j;s),a}$}.

We use this interpretation and the following intuitive idea: when agent $i$ is in site $s$ in state $x_s$ and we decide to send it to site $a$ in state $x_a$, in some sense we are increasing the values of $\rho^i_{(x_s;s),a}$ and $\rho^i_{(x_a;s),a}$, which are the long-term probabilities of such transitions. In particular, we would like to keep the quantities $\bar \rho^i_{(x_j;s),a}$ found to be $0$ in the relaxation as close to $0$ as possible in the final solution. By complementary slackness it is only for these variables that we might have $\bar \gamma^i_{(x_j;s),a}>0$. Hence, when the system is in state $(\mathbf{x};\mathbf{s})$, we associate to each action $\mathbf{a}$ an \emph{index of undesirability}
\begin{align} 	
I((\mathbf{x};\mathbf{s}),\mathbf{a})&=\sum_{\{i \in [N]: s_i \neq a_i\}}^N (\bar \gamma^i_{(x_{s_i};s_i),a_i} + \bar \gamma^i_{(x_{a_i};s_i),a_i}) + \sum_{\{i \in [N]: s_i=a_i\}} \bar \gamma^i_{(x_{s_i};s_i),s_i}, \label{eq: undesirability index}
\end{align}
that is, we sum the reduced costs for the $N$ different projects. 
Then we select an action $\mathbf{a_{pd}} \in \Pi_{[N]}$ that minimizes these indices:
\begin{equation} \label{PD index}
\mathbf{a_{pd}}(\mathbf{x};\mathbf{s}) \in \text{argmin}_{\mathbf{a}} \{I((\mathbf{x};\mathbf{s}),\mathbf{a}) \}.
\end{equation}

We now show that this policy is in fact equivalent to the one-step lookahead policy described earlier. Using the expression for the reduced costs from paragraph \ref{section: relaxation dual}, we can rewrite the indices in (\ref{eq: undesirability index}) more explicitely. The term $\bar \gamma^i_{(x_{s_i};s_i),a_i} + \bar \gamma^i_{(x_{a_i};s_i),a_i}$  of the sum (\ref{eq: undesirability index})  is equal to
\begin{align}
& \bar \lambda^{i}_{s_i, x_{s_i}} + \bar \kappa_{s_i, x_{s_i}} - \sum_{i'=1}^N \bar \xi^{i'}_{s_i} + \sum_{s'=1}^N \bar \xi^{i}_{s'} - \sum_{i'=1}^N \bar \zeta^{i'}_{a_i} + \sum_{a'=1}^N \bar \zeta^{i}_{a'} - \bar \kappa_{a_i, x_{a_i}} \label{eq: invariance} \\ 
& - \alpha \sum_{\tilde x_{a_i}} p^{{1}\{i \leq M\}}_{x_{a_i} \tilde x_{a_i}} \bar \lambda^{i}_{a_i, \tilde x_{a_i}} - r^{{1}\{i \leq M\}}_{a_i} (x_{a_i}) + c_{s_i a_i} {1} \{ i \leq M \}, \nonumber
\end{align}
after cancellation of $\bar \mu^{i}_{s_i, a_i}$, and adding and subtracting $\bar \zeta^{i}_{a_i}$ and $\bar \xi^{i}_{s_i}$. This expression is valid for the terms $\gamma^i_{(x_{s_i};s_i),s_i}$ as well. Now after summation over $i \in [N]$, the first two lines in expression (\ref{eq: invariance}) do not play any role in the minimization (\ref{PD index}). This is obvious for the terms that do not depend on $\mathbf{a}$. For the terms involving the $\zeta^i_{a_j}$, we can write
\[
\sum_{i=1}^N \sum_{i'=1}^N \bar \zeta^{i'}_{a_i} = \sum_{i'=1}^N \sum_{i=1}^N \bar \zeta^{i'}_{a_i} = \sum_{i'=1}^N \sum_{a=1}^N \bar \zeta^{i'}_{a},
\]
the last equality being true since $\mathbf{a}$ is just a permutation of $\{1,\ldots,N\}$. Hence the sums involving the $\bar \zeta^{i}_{a_i}$ cancel (in fact we even see that each individual sum is independent of the choice of $\mathbf{a}$). As for the term $\sum_{i=1}^N \bar \kappa_{a_i, x_{a_i}}$, it is equal to $\sum_{j=1}^N \bar \kappa_{j, x_{j}}$ and so it is independent of the choice of $\mathbf{a} \in \Pi_{[N]}$. We are left with the following optimization problem:
\begin{align*} 
&\mathbf{a_{pd}}(\mathbf{x};\mathbf{s}) \in \text{argmin}_{\mathbf{a}} \left \{ - \sum_{i=1}^N \left( r^{{1}\{i \leq M\}}_{a_i} (x_{a_i}) \right. \right. \left. \left. - c_{s_i a_i} {1} \{ i \leq M \} + \alpha \sum_{\tilde x_{a_i}} p^{{1}\{i \leq M\}}_{x_{a_i} \tilde x_{a_i}} \bar \lambda^{i}_{a_i, \tilde x_{a_i}} \right) \right\},
\end{align*}
which after a sign change is seen to be exactly (\ref{eq: linear assignment}). We have shown the following

\begin{thm} \label{thm: equivalence}
The primal-dual heuristic (\ref{PD index}), based on the interpretation of the reduced costs of the LP relaxation, is equivalent to the one-step lookahead policy (\ref{eq: linear assignment}) assuming the separable approximation (\ref{approximate value}) for the reward-to-go.
\end{thm}

In view of this result, we obtain an alternative way to compute the one-step lookahead policy. The minimization (\ref{PD index}) is again a linear assignment problem. If we can store the $\mathcal{O}(N^3 \, \max_i(|S_i|))$ optimal reduced costs instead of the $\mathcal{O}(N^2 \, \max_i(|S_i|))$ optimal dual variables, there is just a linear cost involved in computing the indices $I((\mathbf{x};\mathbf{s}),\mathbf{a})$ of the problem at each period resulting in an overall $\mathcal{O}(N^3)$ computational cost for the on-line maximization at each period. 

\section{Numerical Experiments} \label{experiments}

Table \ref{experiments table} presents numerical experiments on problems whose characteristics differently affect the performance of the heuristic described in section~\ref{algorithms}. Linear programs are implemented in AMPL and solved using CPLEX. Due to the size of the state space, the expected discounted reward of the heuristics is computed using Monte-Carlo simulations. The computation of each trajectory is terminated after a sufficiently large, but finite horizon: in our case, when $\alpha^t$ times the maximal absolute value of any immediate reward becomes less than $10^{-6}$. To reduce the amount of computation in the evaluation of the policies, we assume that the distribution of the initial states of the sites is deterministic. 

In a given problem, the number $|S_i|$ of states is chosen to be the same for all projects. $c/r$ is the ratio of the average switching cost divided by the average active reward. This is intended to give an idea of the importance of the switching costs in the particular experiment. The switching costs are always taken to be nonnegative. $Z^*$ is the optimal value of the problem, computed using (\ref{RBSC Multi-agent}), when possible. $Z_r$ is the optimal value of the relaxation and so provides an upper bound on the achievable performance. $Z_{osl}$ is the estimated expected value of the one-step lookahead policy. $Z_{g}$ is the estimated expected value of the greedy policy which is obtained by fixing the value of the $\lambda^i_{a_i,x'_{a_i}}$ in (\ref{eq: linear assignment}) to zero, i.e., approximating the reward-to-go by zero. This greedy policy is actually optimal for the MABP with deteriorating active rewards, i.e., such that projects become less profitable as they are worked on \cite[vol. 2, p.69]{bertsekas-dynamicProgramming}. Problem 2 is of this type and shows that the one-step lookahead policy does not perform optimally in general.

Problem 1 is a MABP. The heuristic is not optimal, so we see that we do not recover Gittins' policy. Hence the heuristic is also different from Whittle's in general, which reduces to Gittins' in the MAB case. In problem 3, we add transition costs to problem 2. The greedy policy is not optimal any more, and the one-step lookahead policy performs better in this case. Problem 4 is designed to make the greedy policy underperform: two remote sites have slightly larger initial rewards (taking into account the cost for reaching them), but the active rewards at these sites are rapidly decreasing and the agents are overall better off avoiding these sites. The greedy policy does not take into account the future transition costs incurred when leaving these sites. In this case, it turns out that the one-step lookahead is quasi-optimal.  Problem 7 and 8 are larger scale problems, with up to 30 sites. The relaxation is computed in about 20 minutes on a standard desktop, showing the feasibility of the approach for this range of parameters.

\begin{table*}[!t]		
\caption{Numerical Experiments}
\label{experiments table}
\begin{center}
\begin{tabular*}{0.9\textwidth}{|@{\extracolsep{\fill}}c c c c c c|}		
\hline
Problem 							&	$\alpha$ 	& $Z^*$ 		& $Z_r$ 		& $Z_g$ 	& $Z_{osl}$ 	\\ 
$(N, M, |S_i|, c/r)$ 	&						&						&						&					&							\\ 
\hline
Problem 1 						& 	0.5 		& 	84.69		& 	85.21	 	& 	84.5	& 	84.3			\\ 
(4,1,3,0)							&		0.9			&		299.6 	& 	301.4		&		276 	&		294				\\ 
											&		0.99		&		2614.2	& 	2614		&		2324 	&		2611			\\ 
\hline
Problem 2 						& 	0.5 		& 	84.13		& 	85.14 	& 	84.1 	& 	84.1 			\\ 
(4,1,3,0)							&		0.9			&		231.0		& 	245.1		&		231		&		228				\\ 
											&		0.99		&		1337		& 	1339		&		1337 	&		1336			\\ 
\hline
Problem 3 						& 	0.5 		& 	57.54		& 	59.32 	& 	56.0	&	 	57.3			\\ 
(4,1,3,0.6) 					&		0.9			&		184.5		& 	185.0		&		177		&		183				\\ 
											&		0.99		&		1279		& 	1280		&		1273 	&		1277			\\ 
\hline
Problem 4 						& 	0.5			& 	   			& 	165.7 	& 	115 	& 	165				\\ 
(4,2,5,1.39)						&		0.9			&						& 	767.2		&		661		&		767				\\ 
											&		0.95		&						& 	1518		&		1403 	&		1518			\\ 
\hline
Problem 5 						& 	0.5 		& 	  			& 	39.25 	& 	38.5	& 	36.5 			\\ 
(6,2,4,0)							&		0.9			&						& 	214.0		&		205		&		198				\\ 
											&		0.95		&						& 	431.6		&		414 	&		396				\\ 
\hline
Problem 6 						& 	0.5			& 	  			& 	9.727  	& 	6.93 	& 	8.24 			\\ 
(6,2,4,1.51)						&	  0.9			&						&		62.80 	&		38.0	&		47.0			\\ 
											&		0.95		&						& 	128.7		&		78.0 	&		99.0			\\ 
\hline
Problem 7 						& 	0.5			& 	  			& 	196.5	 	& 	189	 	& 	 194			\\ 
(20,15,3,1.16)				&		0.9			&						&		952.7		&		877		&		 900			\\ 
											&		0.95		&						& 	1899		&		1747 	&		 1776			\\ 
\hline
Problem 8 						& 	0.5			& 	  			& 	 589.4 	& 	566 	& 	 564			\\ 
(30,15,2,2.18)					&	  0.9			&						&		 2833   &	 	2640	&	 	 2641 		\\ 
											&		0.95		&						& 	 5642		&		5218 	&		 5246			\\ 
\hline
\end{tabular*}
\end{center}
\end{table*}

\section{A ``Performance'' Bound} \label{bound}

In this section, we present a result that offers some insight into why we could expect the one-step lookahead policy to perform well if the linear programming relaxation of the original problem is sufficiently tight. We begin with the following


\begin{lem} \label{lem: super-harmonic}
The approximate reward-to-go (\ref{approximate value}) is a feasible solution for the original dual linear program (\ref{eq: dual exact RBSC}).
\end{lem}

\begin{proof}
Consider one constraint in the original dual LP (\ref{eq: dual exact RBSC}), for a fixed state-action tuple $(\mathbf{x},\mathbf{s},\mathbf{a})$. We consider a situation where $s_i \neq a_i$ for all $i \in \{1,\ldots, N\}$. Summing the constraints (\ref{eq: constraint xs,s,a}) over $i$ for the given values of $x_{s_i},s_i,a_i$, we get
\begin{align*}
& \sum_{i=1}^N \lambda^{i}_{s_i, x_{s_i}} + \sum_{i=1}^N \mu^{i}_{s_i, a_i} + \sum_{i=1}^N \kappa_{s_i, x_{s_i}} - \sum_{i=1}^N \sum_{i'=1}^N \zeta^{i'}_{a_i} + \sum_{i=1}^N \sum_{a'=1}^N \zeta^{i}_{a'} \nonumber \\
& = \sum_{i=1}^N \lambda^{i}_{s_i, x_{s_i}} + \sum_{i=1}^N \mu^{i}_{s_i, a_i} + \sum_{i=1}^N \kappa_{s_i, x_{s_i}} \geq 0. 
\end{align*}
The cancellation follows from the discussion preceding theorem \ref{thm: equivalence}. Now summing the constraints (\ref{eq: constraint xa,s,a}) over $i$, we also get 
\begin{align*}
& - \alpha \sum_{i=1}^N \sum_{\tilde x_{a_i}} p^{{1}\{i \leq M\}}_{x_{a_i} \tilde x_{a_i}} \lambda^{i}_{a_i, \tilde x_{a_i}} - \sum_{i=1}^N \mu^{i}_{s_i, a_i} - \sum_{i=1}^N \kappa_{a_i, x_{a_i}} \geq \sum_{i=1}^N r^{{1}\{i \leq M\}}_{a_i} (x_{a_i}) - c_{s_i a_i} {{1}\{i \leq M\}}.
\end{align*}
Finally, we add these two inequalities. We obtain
\begin{align*}
& \sum_{i=1}^N \lambda^{i}_{s_i, x_{s_i}}  - \alpha \sum_{i=1}^N \sum_{\tilde x_{a_i}} p^{{1}\{i \leq M\}}_{x_{a_i} \tilde x_{a_i}} \lambda^{i}_{a_i, \tilde x_{a_i}} \geq \sum_{i=1}^N \left( r^{{1}\{i \leq M\}}_{a_i} (x_{a_i}) - c_{s_i a_i} {{1}\{i \leq M\}} \right),
\end{align*}
which is the inequality obtained by using the vector (\ref{approximate value}) in the constraints of (\ref{eq: dual exact RBSC}). 

The case where $s_i=a_i$ for some $i$ is almost identical, considering the constraints (\ref{eq: constraint xs,s,s}) for the corresponding indices.
\end{proof}

In the following theorem, the occupation measure $F_\alpha(\nu,\tilde{u})$ is a vector of size $|\mathcal{S}|$, representing the discounted infinite horizon frequencies of the states under policy $\tilde u$ and initial distribution $\nu$ \cite{derman-LP}. The proof of the theorem follows from the analysis presented in \cite{deFarias-ADP}, see \cite{leny06_RBSC_CDC} for more details.


\begin{thm}
Let $\nu$ be an initial distribution on the states, of the product form (\ref{eq: initial distribution - product}). Let $J^*$ be the optimal reward function, $\tilde J$ be an approximation of this reward function which is feasible for the LP (\ref{eq: dual exact RBSC}), and $\tilde{u}$ be the associated one-step lookahead policy. Let $F_\alpha(\nu,\tilde{u})$ and $J_{\tilde{u}}$ be the occupation measure vector and the expected reward associated to the policy $\tilde{u}$. Then

\begin{align} \label{upperBound}
\nu^T(J^*-J_{\tilde{u}}) \leq \frac{1}{1-\alpha} F_\alpha(\nu,\tilde{u})^T (\tilde{J}-J^*).
\end{align}
\end{thm}

From lemma \ref{lem: super-harmonic}, the theorem is true in particular for $\tilde J$ formed according to (\ref{approximate value}). In words, it says that starting with a distribution $\nu$ over the states, the difference in expected rewards between the optimal policy and the one-step lookahead policy is bounded by a weighted $l^1$-distance between the estimate $\tilde{J}$ used in the design of the policy and the optimal value function $J^*$. The weights are given by the occupation measure of the one-step lookahead policy. It provides some motivation to obtain a good approximation $\tilde J$, i.e., a tight relaxation, which was an important element of this paper.

\section{Acknoledgements}
The authors thank Michel Goemans, Sean Meyn and John Tsitsiklis for valuable comments. This work was supported by Air Force - DARPA - MURI award 009628-001-03-132 and Navy ONR award N00014-03-1-0171. 



\bibliographystyle{elsart-num-sort}
\bibliography{bandits}

\begin{thebibliography}{10}
\expandafter\ifx\csname url\endcsname\relax
  \def\url#1{\texttt{#1}}\fi
\expandafter\ifx\csname urlprefix\endcsname\relax\def\urlprefix{URL }\fi

\bibitem{asawa-penalties}
M.~Asawa, D.~Teneketzis, Multi-armed bandits with switching penalties, IEEE
  transactions on automatic control 41~(3) (1996) 328--348.

\bibitem{bertsekas-dynamicProgramming}
D.~Bertsekas, Dynamic Programming and Optimal Control, vol. 1 and 2, 2nd ed.,
  Athena Scientific, 2001.

\bibitem{bertsimas-RB}
D.~Bertsimas, J.~Ni{\~n}o-Mora, Restless bandits, linear programming
  relaxations, and a primal-dual index heuristic, Operations Research 48 (2000)
  80--90.

\bibitem{derman-LP}
C.~Derman, Finite State Markovian Decision Processes, Academic Press, 1970.

\bibitem{gittinsJones}
J.~Gittins, D.~Jones, A dynamic allocation index for the sequential design of
  experiments, in: J.~Gani (ed.), Progress in Statistics, North-Holland,
  Amsterdam, 1974, pp. 241--266.

\bibitem{Glazebrook06_machine_index}
K.~D. Glazebrook, D.~Ruiz-Hernandez, C.~Kirkbride, Some indexable families of
  restless bandit problems, Advances in Applied Probability 38~(3) (2006)
  643--672.

\bibitem{Guha07_RBApprox}
S.~Guha, K.~Munagala, P.~Shi, On index policies for restless bandit problems
  (2007).
\newline\urlprefix\url{http://arxiv.org/abs/0711.3861}

\bibitem{jun-surveypaper}
T.~Jun, A survey on the bandit problem with switching costs, De Economist
  152~(4) (2004) 513--541.

\bibitem{jun05-investments}
T.~Jun, Costly switching and investment volatility, Economic Theory 25~(2)
  (2005) 317--332.

\bibitem{koole-server}
G.~Koole, Assigning a single server to inhomogeneous queues with switching
  costs, Theoretical Computer Science 182 (1997) 203--216.

\bibitem{leny06_RBSC_CDC}
J.~{Le Ny}, M.~Dahleh, E.~Feron, Multi-agent task assignment in the bandit
  framework, in: Proccedings of the 45th IEEE Conference on Decision and
  Control, San Diego, CA, 2006.

\bibitem{leny06_RBSC_ACC}
J.~{Le Ny}, E.~Feron, Restless bandits with switching costs: linear programming
  relaxations, performance bounds and limited lookahead policies, in:
  Proceedings of the American Control Conference, Minneapolis, MN, 2006.

\bibitem{ninoMora-MABPSC}
J.~Ni{\~n}o-Mora, A faster index algorithm and a computational study for
  bandits with switching costs, INFORMS Journal on Computing. Accepted.

\bibitem{papadimitriou-complexityQueuing}
C.~Papadimitriou, J.~Tsitsiklis, The complexity of optimal queueing network
  control, Mathematics of Operations Research 24~(2) (1999) 293--305.

\bibitem{deFarias-ADP}
D.~{Pucci de Farias}, B.~{Van Roy}, The linear programming approach to
  approximate dynamic programming, Operations Research 51~(6) (2003) 850--856.

\bibitem{puterman-MDP}
M.~Puterman, Markov Decision Processes, Wiley, 1994.

\bibitem{schrijver-book}
A.~Schrijver, Combinatorial Optimization - Polyhedra and Efficiency, Springer,
  2003.

\bibitem{vanOyen-indexPolicies}
M.~{Van Oyen}, D.~Pandelis, D.~Teneketzis, Optimality of index policies for
  stochastic scheduling with switching penalties, Journal of Applied
  Probability 29 (1992) 957--966.

\bibitem{varaiya-extensions}
P.~Varaiya, J.~Walrand, C.~Buyukkoc, Extensions of the multiarmed bandit
  problem: the discounted case, IEEE transactions on automatic control 30~(5)
  (1985) 426--439.

\bibitem{whittle-restless}
P.~Whittle, Restless bandits: activity allocation in a changing world, Journal
  of Applied Probability 25A (1988) 287--298.

\end{thebibliography}







\end{document}